\documentclass[11pt]{amsart}

\usepackage{amsmath,amsfonts,amssymb,amscd}
\begin{document}
\renewcommand{\thefootnote}{\fnsymbol{footnote}}
\pagestyle{plain}

\title{The Donaldson-Futaki invariant for\\
sequences of test configurations}
\author{Toshiki Mabuchi${}^*$}
\maketitle
\footnotetext{ ${}^{*}$Supported 
by JSPS Grant-in-Aid for Scientific Research (A) No. 20244005.}
\abstract
In this note, given a polarized algebraic manifold $(X,L)$,
we define the Donaldson-Futaki invariant $F_1 (\{\mu_i\})$
for a sequence $\{\mu_i\}$ of 
test configurations for $(X,L)$ of exponents $\ell_j$ satisfying
$$
\ell_j \to  \infty, \quad\text{ as $j \to \infty$.}
$$
This then allows us to define a strong version of K-stability or
K-semistability for $(X,L)$.  In particular, $(X,L)$ will be shown to be K-semistable in this strong sense if 
the polarization class $c_1(L)_{\Bbb R}$
admits a constant scalar curvature K\"ahler metric. 
\endabstract

\section{Introduction}

By a {\it polarized algebraic manifold} $(X,L)$, we mean
a pair of a nonsingular irreducible projective variety $X$, 
defined over $\Bbb C$, and a very ample line bundle $L$
over $X$. 
Replacing $L$ by its suitable multiple if necessary, 
we may assume 
$$
H^q(X, L^{\otimes \ell }) = \{0\}, \qquad q= 1,2,\dots,n,
$$
%
%
for all positive integers $\ell$, where $n:= \dim_{\Bbb C} X$.
In this note, we fix once for all such a pair $(X,L)$. For the affine line 
 $ \Bbb A^1 := \{z\in \Bbb C\}$, we consider the algebraic torus $T := \Bbb C^*$ acting on $\Bbb A^1$ by multiplication of complex numbers,
$$
T\times \Bbb A^1 \to \Bbb A^1,
\qquad (t, z) \mapsto
t z.
\leqno{(1.1)}
$$
For each positive integer $\ell$, 
the space $V_{\ell}:= H^0(X,L^{\otimes \ell})$ 
is assumed to have a Hermitian metric $\rho_{\ell}$ expressible in the form (1.6) below. Let
 $$
 \Phi_{\ell} \,:\, X \, \hookrightarrow \, \Bbb P^*(V_{\ell})
 $$
be the Kodaira embedding of $X$ associated to the complete linear system $|L^{\otimes \ell}|$ on $X$.
By fixing an $\ell$, we consider an algebraic group homomorphism
$$
\psi \, : \,T \,\to \,\operatorname{GL}(V_{\ell})
$$
such that the maximal compact subgroup $S^1 \subset \Bbb C^* \,(= T)$
acts isometrically on $(V_{\ell}, \rho_{\ell})$.
Let $\mathcal{X}^{\psi}$ be the irreducible algebraic subvariety of $\Bbb A^1 \times \Bbb P^* (V_{\ell})$ obtained as 
the closure of $\cup_{z\in T} \mathcal{X}^{\psi}_z$ in $\Bbb A^1 \times \Bbb P^* (V_{\ell})$ 
by setting
$$
\mathcal{X}^{\psi}_z := \{z\}\times\psi (z)  \Phi_{\ell} (X),
\qquad z \in T,
$$
where the element $\psi (z)  $ in $\operatorname{GL}(V_{\ell})$ 
acts naturally on  the 
set $\Bbb P^* (V_{\ell})$ of all hyperplanes in $V_{\ell}$ passing through the origin.
 We then consider the map
$$
\pi : \mathcal{X}^{\psi} \to \Bbb A^1
$$ 
induced by the projection of $\Bbb A^1 \times \Bbb P^* (V_{\ell})$ to the first factor $\Bbb A^1$.
Moreover, for the hyperplane bundle $\mathcal{O}_{\Bbb P^*(V_{\ell})}(1)$ on $\Bbb P^*(V_{\ell})$, we put 
$$
\mathcal{L}^{\psi}\, :=\,\operatorname{pr}_2^*\mathcal{O}_{\Bbb P^*(V_{\ell})}(1)_{|\mathcal{X}^{\psi}},
$$
where $\operatorname{pr}_2 : \Bbb A^1 \times \Bbb P^* (V_{\ell}) \to \Bbb P^* (V_{\ell})$
is the projection to the second factor.
For the dual vector space $V_{\ell}^*$ of $V_{\ell}$,
the $\Bbb C^*$-action on $\Bbb A^1 \times V_{\ell}^*$ defined by
$$
\Bbb C^* \times (\Bbb A^1 \times V_{\ell}^*)\to \Bbb A^1 \times V_{\ell}^*,
\quad (t, (z, p))\mapsto  (tz, \psi (t) p),
$$
naturally induces a $\Bbb C^*$-action on $\Bbb A^1 \times \Bbb P^*(V_{\ell})$ and $\mathcal{O}_{\Bbb P^*(V_{\ell})}(-1)$, where $\operatorname{GL}(V_{\ell})$ 
acts  on $V_{\ell}^*$ 
by contragradient representation.
 This then induces $\Bbb C^*$-actions on $\mathcal{X}^{\psi}$ and $\mathcal{L}^{\psi}$, and $\pi : \mathcal{X}^{\psi} \to \Bbb A^1$ is a projective morphism with relative very ample line bundle 
$\mathcal{L}^{\psi}$ such that
$$
(\mathcal{X}^{\psi}_z, \mathcal{L}_z^{\psi})\; \cong \; (X,L^{\otimes \ell}),
\qquad z \neq 0,
$$
where for each $z \in \Bbb A^1$, we denote by $\mathcal{L}_z^{\psi}$ the restriction of $\mathcal{L}^{\psi}$ to 
the scheme-theoretic fiber $\mathcal{X}^{\psi}_z$ of $\mathcal{X}^{\psi}$ over $z$.
Then a  triple $({\mathcal{X}}, {\mathcal{L}}, \psi )$ is called 
a {\it 
test configuration for $(X,L)$}, if we have both 
$$
\mathcal{X} = \mathcal{X}^{\psi}\quad\text{ and  }\quad\mathcal{L}= \mathcal{L}^{\psi}
$$
for some $\psi$ as above,
where 
$\ell$  is called 
the {\it exponent} of the test configuration $({\mathcal{X}}, {\mathcal{L}}, \psi )$.
We consider the projective linear representation 
$$
\psi^{\operatorname{PGL}}: \,T \,\to \,\operatorname{PGL}(V_{\ell})
$$ 
induced by $\psi$. 
Then a 
test configuration $({\mathcal{X}}, {\mathcal{L}}, \psi )$
for $(X,L)$ is called {\it trivial}, if 
$\psi^{\operatorname{PGL}}$ 
is a trivial homomorphism.
We now consider 
the set $\mathcal{M}$ of all sequences $\{\mu_j \}$ of 
test configurations 
$$
\mu_j  \, =\, (\mathcal{X}_j, \mathcal{L}_j, \psi_j ), \qquad j =1,2,\dots,
\leqno{(1.2)}
$$
for $(X,L)$  
such that the exponent $\ell_j$ of the 
test configuration $\mu_j$ satisfies the following 
growth condition:
$$
\text{$\ell_j \to +\infty$, \; as $j \to \infty$.}
\leqno{(1.3)}
$$

\smallskip\noindent
Assume now that $\ell_1 =1$ for simplicity.
In Section 2, for every element $\{\mu_j\}$ 
in $\mathcal{M}$, we shall define
the Donaldson-Futaki invariant 
$$
F_1 (\{\mu_j\})\; \in \; \Bbb R \cup \{-\infty\}.
$$ 
This definition is justified in the sense that $F_1 (\{\mu_j\}) \in \Bbb R\cup \{-\infty\}$ can be viewed as the Donaldson-Futaki invariant $F_1 (\mu_{\infty})$ for the limit 
$\mu_{\infty}$ 
of a suitable subsequence of $\{\mu_j\}$
in some completion of the moduli space of 
test configurations for $(X,L)$. 
This will be shown in \cite{MN} in general, while the most important case is the following:

\medskip\noindent
For a $T$-action on an irreducible algebraic variety $\mathcal{X}$,   
we assume that there exists a $T$-equivariant projective morphism 
$$
\pi : \mathcal{X} \to \Bbb A^1
$$ 
with a relatively very ample line bundle $\mathcal{L}$ on  the fiber space $\mathcal{X}$ over $\Bbb A^1$ such that the $T$-action on $\mathcal{X}$ 
lifts to a $T$-linearization of $\mathcal{L}$. Let $(\mathcal{X},\mathcal{L})$ be a {\it test configuration 
for \/$(X,L)$ of exponent $1$ in Donaldson's sense} (see \cite{D1}), i.e., we have an algebraic isomorphism
$$
(\mathcal{X}_z,\mathcal{L}_z )\;  \cong \; (X, L), 
\qquad 0 \neq  z \in \Bbb A^1,
$$
where $\mathcal{L}_z$ denotes the restriction of $\mathcal{L}$ to $\mathcal{X}_z : = \pi^{-1}(z)$.
Then by the affirmative solution of equivariant Serre's conjecture, the direct image sheaves
$$
E^{j} \, := \; \pi_*\mathcal{L}^{\otimes \ell_j},
\qquad j = 1,2,\dots,
$$
can be viewed as the trivial vector bundles over $\Bbb A^1$ in terms of a $T$-equivariant isomorphism 
of vector bundles 
$$
E^{j} \cong \Bbb A^1\times E^{j}_0,
\leqno{(1.4)}
$$
where the Hermitian metric $\rho_{\ell_j}$ on 
the fiber $E^{j}_1 \; (= V_{\ell_j})$ over $1$
is chosen to be $\rho_h^{(\ell_j)}$  in (1.6)  inducing a Hermitian metric on the central 
fiber $E^{\ell}_0$ which is preserved by the action of $S^1\subset T$ on $E^{j}_0$ (see \cite{D2}).
Here the $T$-action on $\Bbb A^1\times E^{j}_0$ is induced by the $T$-action (1.1) on $\Bbb A^1$
and a natural $T$-representation
$$
\psi_j \,: \; T \to \operatorname{GL}(V_{j}) \; (= \operatorname{GL}(E^{j}_0)),
$$
under the identification of $E^{j}_0$ with $E^{j}_1 \; (= V_{j})$ by the isomorphism (1.4).
In view of $\ell_1 = 1$, we write $\psi_1$ simply as $\psi$. We now assume that 
the test configuration $(\mathcal{X},\mathcal{L})$ is {\it nontrivial}, i.e., 
the homomorphism $\psi^{\operatorname{PGL}}$ is nontrivial.
Put $\lambda = c_1 (X)^n[X]/\pi $.
The following will be proved in a separate paper:

\medskip\noindent
{\bf  Fact} (cf. \cite{M8}).
{\em For the test configuration $(\mathcal{X},\mathcal{L})$ as above for $(X,L)$
in Donaldson's sense, consider the sequence of 
test configurations
$\mu_j = (\mathcal{X}_j, \mathcal{L}_j, \psi_j)$ of exponents $\ell_j$, $j = 1,2,\dots,$ defined by
$$
\mathcal{X}_j := \mathcal{X}\;\;\text{ and }\;\;\mathcal{L}_j := \mathcal{L}^{\otimes \ell_j},
$$
where $\ell_j$, $j =1,2,\dots$, are positive integers satisfying $\ell_1 =1$ and $(1.3)$. 
Replacing 
the sequence $\{\ell_j\}$ by its suitable subsequence if necessary, in view of  \/$[5;$ Remark $4.6\,]$
and  \/$[6;$ Lemma $2.6\,]$,
we may assume that
$|b_j|_1^{}$ in \/$(2.2)$ converges to a nonnegative real number $\beta$.
If \/$\beta > 0$, then for the sequence $\{\mu_j \}$ 
of 
test configurations, we have
$$
F_1 (\{\mu_j\}) \; =\; \beta^{-1}\,\lambda\,F_1 (\mathcal{X},\mathcal{L} ),
\leqno{(1.5)}
$$
where $F_1 (\mathcal{X},\mathcal{L} )$ 
on the right-hand side  is the Donaldson-Futaki invariant of the test configuration $(\mathcal{X},\mathcal{L} )$ in Donaldson's sense.
} 

\medskip
Now we come back to a general situation.
Fix a Hermitian metric $h$ for $L$ such that $\omega := c_1(L;h)$ is K\"ahler.
For each positive integer $\ell$,
we choose $\rho_{\ell}$ as the Hermitian metric $\rho_h^{(\ell )}$ on $V_{\ell}$ defined by
$$
\rho_h^{(\ell )}(\sigma', \sigma'' ) \;:=\; \int_X (\sigma', \sigma'')^{}_h\, \omega^n,
\qquad \sigma',\sigma'' \in V_{\ell},
\leqno{(1.6)}
$$
where $(\sigma',\sigma'' )_h$ denotes the pointwise Hermitian inner product of 
$\sigma$ and $\sigma'$ by the $\ell$-multiple of $h$.
We shall now define a strong version of 
K-stability which will be studied in detail in a forthcoming paper \cite{M6}.

\medskip\noindent
{\em Definition\/ $1.7$}. (1) 
A polarized algebraic manifold $(X,L)$ 
 is called 
{\it strongly K-semistable}, 
if $F_1(\{\mu_j\} ) \leq 0$ for all 
$\{\mu_j\} \in 
\mathcal{M}$.

\medskip\noindent
(2) A strongly K-semistable polarized algebraic manifold $(X,L)$ 
is called 
{\it strongly $K$-stable}, 
if  for every $\{\mu_j\} \in \mathcal{M}$ satisfying $ F_1(\{\mu_j\} ) = 0$,
there exists a $j_0$ such that
the 
test configurations $\mu_j$ are trivial for all $j\geq j_0$.


 
\medskip
We here observe that these definitions of  strong K-semistability and strong K-stability are independent of the choice 
of $h$ above (cf.~\cite{MN}).
Finally in Section 4, we shall prove the following:

\medskip\noindent
{\bf Theorem.} 
{\em If the polarization class $c_1(L)_{\Bbb R}$ admits a constant scalar curvature K\"ahler metric, then 
$(X,L)$ is strongly K-semistable.} 

\medskip
In \cite{M6}, we shall show a stronger result that
the polarized algebraic manifold $(X,L)$ in the Theorem above is actually strongly K-stable.

\section{The Donaldson-Futaki invariant $F_1$ on $\mathcal{M}$}

Let $\mu_j = (\mathcal{X}_j, \mathcal{L}_j, \psi_j)$, $j = 1,2,\dots$,  be a sequence of test configurations 
for $(X,L)$ as in (1.2), and let $\ell_j$ be the exponent of $\mu_j$.
By abuse of terminology, we write the vector space $V_{\ell_j} = H^0(X, L^{\otimes \ell_j} )$ simply as $V_j$.
Put $N_j := \dim V_j$.
Then by setting $d_j := \ell_j^{\,n} c_1(L)^n [X]$, we define
$$
W_{j} \; := \; \{\operatorname{Sym}^{d_j }(V_j )\}^{\otimes n+1},
$$
where $\operatorname{Sym}^{d_j}(V_j )$ denotes the $d_j$-th symmetric tensor product  of $V_j$.
The dual space $W_{j}^*$ of $W_{j}$ admits the Chow norm (cf.~\cite{Zh})
$$
W_{j}^* \owns w \;\mapsto \;\| w\|^{}_{\operatorname{CH}(\rho^{}_{\ell_j})} \in \Bbb R^{}_{\geq 0},
$$
associated to the Hermitian metric $\rho^{}_{\ell_j}$ on $V_j$ as in (1.6). For the Kodaira embedding
$\Phi_{j}: X \hookrightarrow  \Bbb P^*(V_{j})$ associated to the complete linear 
system $|L^{\otimes \ell_j}|$ on $X$, we consider the Chow form
$$
0 \neq \hat{X}_{j}\in W_{j}^*
$$ 
for the irreducible reduced algebraic cycle 
$\Phi_{j}(X)$ on $\Bbb P^*(V_{j})$ such that the corresponding point $[\hat{X}_{j}]$ in $\Bbb P^*(W_j )$
is the Chow point for the cycle $\Phi_{j}(X)$.
For the homomorphism
 $\psi_{j}: T \to \operatorname{GL}(V_j)$,
taking the real Lie subgroup 
$$
T_{\Bbb R}  = \{t \in T\,;\, t\in \Bbb R_+\}
$$
of the algebraic torus $T = \{t \in \Bbb C^*\}$,
we now define a Lie group homomorphism 
$\psi_{j}^{\operatorname{SL}}: {T}_{\Bbb R} \to \operatorname{SL}(V_j)$ by
$$
\psi_{j}^{\operatorname{SL}}(t) \; :=\; \frac{\psi_{j}(t)}{\det (\psi_{j}(t))^{1/N_j}},
\qquad t\in {T}_{\Bbb R}.
$$
Then for a suitable orthonormal basis $\{\sigma_{\alpha}\,;\, \alpha =1,2,\dots,N_{j}\}$ for $V_{j}$,
there exist rational numbers $b_{j,\alpha}$, $\alpha =1,2, \dots,N_{j}$ such that
the homomorphism $\psi_{j}^{\operatorname{SL}}: {T}_{\Bbb R} \to \operatorname{SL}(V_j)$ 
is written in the form
$$
\psi_{j}^{\operatorname{SL}}({t})\cdot \sigma_{\alpha} \; =\; {t}_{}^{-b^{}_{j,\alpha}}\sigma_{\alpha},
\qquad \alpha =1,2,\dots,N_{j},  
\leqno{(2.1)}
$$ 
where ${t}\in \Bbb R_+$ is arbitrary. 
Then for $b_j := (b^{}_{j,1}, b^{}_{j,2}, \dots, b^{}_{j,N_{j}}) \in \Bbb R^{N_j}$ 
as above, 
we shall now define its norms
$$
\begin{cases}
&|b_j |_1\;  := \,
\Sigma_{\alpha =1}^{N_{j}}\, |b^{}_{j,\alpha}|/\ell^{n+1}_j,\\
& |b_j |_{\infty}\;  := \;\max \{ |b_{j,1}|, |b_{j,2}|, \dots, |b_{j,N_{j}}|\}/\ell_j, 
\end{cases}
\leqno{(2.2)}
$$
where
by the definition of $b^{}_{j}$ we have the equality
$$
b^{}_{j,1}+b^{}_{j,2}+\dots +b^{}_{j,N_j}\; =\; 0.
\leqno{(2.3)}
$$
Since $\operatorname{SL}(V_j)$ 
acts naturally on $W_{j}^*$, writing each $t\in T_{\Bbb R}$ as $t = \exp (s/|b_j|_{\infty})$ for some $s\in \Bbb R$ in the case $|b_j|_{\infty}\neq 0$, we
define a real function $f_{ j}$ on $\Bbb R$ by
$$
f_{j}(s) \; :=\; 
\begin{cases}
& (|b_j|^{}_1\ell_j^n)^{-1}|b _j|_{\infty}
\log \| \psi_{j}^{\operatorname{SL}}(t)\cdot \hat{X}_{j} \|^{}_{\operatorname{CH}(\rho^{}_{\ell_j})}, 
\quad \,\text{ if $\psi^{\operatorname{PGL}}_j \neq 1$},\\
& \;\ell_j^{-n}\log \|  \hat{X}_{j} \|^{}_{\operatorname{CH}(\rho^{}_{\ell_j})},\quad\quad\; \qquad\qquad\qquad\quad\,
\,\text{ if $\psi^{\operatorname{PGL}}_j = 1$},
\end{cases}
$$
where $s$ runs through the set of all real numbers.
Put $\dot{f}_{j}:= df_{j}/ds$.
Since $h$ is fixed, the derivative $\dot{f}_{j}(0)$ at $s = 0$ is bounded from above by a constant 
$C>0$ independent
of $j$ (see Section 3).
Define $F_{1} (\{\mu_{j}\} )\in \Bbb R \cup \{-\infty \}$ by 
$$
F_{1} (\{\mu_{j}\} )\, :=\; \lim_{s\to -\infty} \{
\varliminf_{j\to \infty}  \dot{f}_{j}(s)\}
\; 
\leq\; C.
\leqno{(2.4)}
$$
This is well-defined, because the function $\varliminf_{j\to \infty}  \dot{f}_{j}(s)$  is non-decreasing  in $s$ by convexity 
of the functions ${f}_{j}$ (cf. \cite{Zh}; see also \cite{M1}, Theorem 4.5).


\section{Boundedness from above for $\dot{f}_{j}(0)$} 

If $\psi^{\operatorname{PGL}}_j = 1$, then $f_j$ is a constant function, so that $\dot{f}_j(0) =0$.
Hence we may assume that $\psi^{\operatorname{PGL}}_j \neq 1$.
For the K\"ahler metric $\omega = c_1(L;h)$ on $X$, let
$B_{j}(\omega )$ denote the $\ell_j$-th asymptotic Bergman kernel defined by
$$
B_{j}(\omega ) \, :=\;\Sigma_{\alpha =1}^{N_j} \,  |\sigma_{\alpha} |_h^{\,2}
\leqno{(3.1)}
$$
for the orthonormal basis $\{\,\sigma_{\alpha}\,;\, \alpha = 1,2,\dots, N_j\,\}$ of $V_j$ 
as in Section 2, 
where we put $|\sigma_{\alpha} |_h^{\,2} := (\sigma_{\alpha}, \sigma_{\alpha})_h$ for simplicity.
We now consider the Kodaira embedding
$\Phi_j : X \hookrightarrow \Bbb P^{N_j-1}(\Bbb C) = \{(z_1:z_2: \dots :z_{N_j})\}$
of $X$ in the projective space $\Bbb P^{N_j-1}(\Bbb C )$ defined by 
$$
\Phi_j (x) := (\sigma_1 (x): \sigma_2 (x): \dots : \sigma_{N_j}(x)), \;\quad x \in X,
$$
where we identify $\Bbb P^*(V_j)=\Bbb P^{N_j-1}(\Bbb C )$.
Let $\omega_{\operatorname{FS}}: = (\sqrt{-1}/2\pi)\partial\bar{\partial}\log \Sigma_{\alpha =1}^{N_j} |z_{\alpha}|^2$  be the Fubini-Study form on $\Bbb P^{N_j-1}(\Bbb C ) $.
 Then by \cite{Ze}, we see that
\begin{align*}
&|B_j(\omega ) \,- \, (1/n!)\,\ell_j^{\,n}  | \; \leq \;  C_1 \ell_j^{\,n-1},\tag{3.2}
\\
&(\ell_j -C_2) \omega \leq \Phi_j^*\omega_{\operatorname{FS}} 
\leq  (\ell_j + C_2) \omega,\tag{3.3}
\end{align*}
where $C_i$, $i=1,2,\dots$, denote positive real constants  independent  of 
the choice of $j$, $\alpha$ and $t\in \Bbb C^*$ throughout this paper. 
By \cite{Zh} (see also \cite{M1}),
$$
\dot{f}_{j}(0)\; =\; (|b_j|^{}_1\ell^n_j )^{-1}\int_X
B_j(\omega )^{-1}
(\Sigma_{\alpha =1}^{N_j}\, b^{}_{j,\alpha} |{\sigma}_{\alpha}|_h^2)\,\Phi_j^*\omega^{\,n}_{\operatorname{FS}},
\leqno{(3.4)}
$$
where the norm $|b_j|_1$ is as in (2.2).
Then 
$(|b_j|^{}_1\ell_j )^{-1}\int_X \Sigma_{\alpha =1}^{N_j}\, b^{}_{j,\alpha} |{\sigma}_{\alpha}|_h^2\omega^n = 0$ 
by (2.3). Hence by setting
$$
\begin{cases}
&R_1:= (|b_j|^{}_1 \ell_j^n )^{-1}\int_X
B_j(\omega )^{-1}
(\Sigma_{\alpha =1}^{N_j}\, b^{}_{j,\alpha} |{\sigma}_{\alpha}|_h^2)\,
\{\Phi_j^*\omega^{\,n}_{\operatorname{FS}} - (\ell_j \omega )^n\},\\
&R_2:= (|b_j|^{}_1\ell_j^n)^{-1}\int_X
\{B_j(\omega )^{-1}- (\ell_j^{\,n}/n!)^{-1}\}
(\Sigma_{\alpha =1}^{N_j}\, b^{}_{j,\alpha} |{\sigma}_{\alpha}|_h^2)\,
(\ell_j \omega )^n,
\end{cases}
$$
we see from (3.4) that $\dot{f}_{j}(0) = R_1+R_2$.
Now by (3.2) and (3.3), 
\begin{align*}
R_1\; 
&
\leq \; \ell_j\,(\Sigma_{\alpha =1}^{N_j}\, |b^{}_{j,\alpha}|)^{-1}\int_X
B_j(\omega )^{-1}
(\Sigma_{\alpha =1}^{N_j}\, |b^{}_{j,\alpha}|\, |{\sigma}_{\alpha}|_h^2)\,
|\Phi_j^*\omega^{\,n}_{\operatorname{FS}} - (\ell_j \omega )^n| \\
& 
\leq \; C_3\, \ell_j\,(\Sigma_{\alpha =1}^{N_j}\, |b^{}_{j,\alpha}|)^{-1}\int_X
(\ell_j^{\,n}/n! )^{-1}
(\Sigma_{\alpha =1}^{N_j}\, |b^{}_{j,\alpha}|\, |\sigma_{\alpha}|_h^2)\,
\ell^{\,n-1}_j \omega^n   \\
& 
=\; n!\, C_3 \,(\Sigma_{\alpha =1}^{N_j}\, |b^{}_{j,\alpha}|)^{-1}\int_X
(\Sigma_{\alpha =1}^{N_j}\, |b^{}_{j,\alpha}|\, |\sigma_{\alpha}|_h^2)\,\omega^n
\; =\; n!\, C_3.
\end{align*}
On the other hand, we see from (3.2) that
\begin{align*}
R_2 \; 
&\leq \; C_4\,\ell_j\,(\Sigma_{\alpha =1}^{N_j}\, |b^{}_{j,\alpha}|)^{-1}\int_X
 \ell_j^{\,-n-1}
(\Sigma_{\alpha =1}^{N_j}\, b^{}_{j,\alpha}\, |{\sigma}_{\alpha}|_h^2)\,
(\ell_j \omega )^n\\
& 
\leq \; C_4 \,(\Sigma_{\alpha =1}^{N_j}\, |b^{}_{j,\alpha}|)^{-1}\int_X
(\Sigma_{\alpha =1}^{N_j}\, |b^{}_{j,\alpha}|\, |\sigma_{\alpha}|_h^2)\,\omega^n
\; =\; C_4.
\end{align*}
Thus we obtain $\dot{f}_{j}(0) = R_1+R_2\leq n!\,C_3  + C_4 = C_5$, 
as required.

\section{Proof of Theorem}
 
In this section, we use the same notation as in Section 3 by choosing
a K\"ahler metric of constant scalar curvature
as our $\omega = c_1(L;h)$.
Associated to this $h$, we take $\rho_h^{(\ell )}$ in (1.6)
as the metric $\rho_{\ell}$ on $V_{\ell}$.
The scalar curvature 
$$
S_{\omega} \; =\;  n\, c_1(M)c_1(L)^{n-1}[X]/c_1(L)^n[X]
$$ 
of $\omega$ is a constant function independent of the choice of $j$. 
Now by \cite{Lu} and \cite{Ze}, the $\ell_j$-th asymptotic Bergman kernel 
$B_j (\omega )$ in (3.1) is written as
$$
B_j (\omega )\; =\; (1/n!)\, \{ \ell_j^n \,+\,(S_{\omega}/2 ) \ell_j^{n-1} \,+\, R^{}_j\,
\ell_j^{n-2}\},
\leqno{(4.1)}
$$
where $R_j$ is a function on $X$ satisfying $\|R_j\|^{}_{C^2(X)}\leq C^{}_6 $.
Now for each integer $r$, let $O(\ell_j^{\,r})$ denote a function $\varphi$ or a real $2$-form $\theta$ on $X$ such that $|\varphi |^{}\leq C^{}_7 \ell_j^{\,r}$ or $-C^{}_7 \ell_j^{\,r}\omega
\leq \theta \leq C^{}_7 \ell_j^{\,r}\omega$, respectively,  for some $C_7$. Then by taking 
$(\sqrt{-1}/2\pi) \partial\bar{\partial}\log$ of both sides of (4.1), we obtain
$$
\Phi_j^*\omega_{\operatorname{FS}} - \ell_j \omega\; =\; O(\ell_j^{-2}).
\leqno{(4.2)}
$$
Hence by (4.1) and (4.2), we see from (3.4) that
\begin{align*}
\dot{f}_{j}(0)\; &=\; (|b_j|_1^{}\ell_j^n )^{-1}\int_X
\frac{n!\,\Sigma_{\alpha =1}^{N_j}\, b^{}_{j,\alpha} |{\sigma}_{\alpha}|_h^2}{   \ell_j^n \,+\,(S_{\omega}/2 ) \ell_j^{n-1} + O(\ell_j^{n-2})}  \,\{\ell_j\omega + O(\ell_j^{-2})\}^{\,n}\\
&=\; \frac{n!\,(|b_j|_1^{}\ell^n_j)^{-1}}{\ell_j^n \,+\,(S_{\omega}/2 ) \ell_j^{n-1} }
\int_X \frac{\Sigma_{\alpha =1}^{N_j}\, b^{}_{j,\alpha} |{\sigma}_{\alpha}|_h^2}{  1 \,+\,O(\ell_j^{-2})} 
\,\{\ell_j\omega + O(\ell_j^{-2})\}^{\,n}.
\end{align*}
Since $\int_X \Sigma_{\alpha =1}^{N_j}\, b^{}_{j,\alpha} |{\sigma}_{\alpha}|_h^2\, \omega^n = 
\Sigma_{\alpha=1}^{N_j} b^{}_{j,\alpha}= 0$ by (2.3), it now follows that
\begin{align*}
|\dot{f}_{j}(0)| \; &\leq \; C_8\,(|b_j|_1^{}\ell_j^n)^{-1}\int_X (\Sigma_{\alpha =1}^{N_j}\, |b^{}_{j,\alpha}|\, |{\sigma}_{\alpha}|_h^2)\,O(\ell_j^{-2})\,\omega^n\\
&=\; C_8\,(|b_j|_1^{}\ell_j^n)^{-1}\,(\Sigma_{\alpha =1}^{N_j}\, |b^{}_{j,\alpha}|)\,O(\ell_j^{-2})
\; =\; O(\ell_j^{-1}),
\end{align*}
where the last equality follows from (2.2). Then by $\lim_{j\to \infty}\ell_j^{-1} =0$,
 we obtain $\lim_{j\to\infty} \dot{f}_{j}(0) 
=0$. 
Since $\varliminf_{j\to \infty}  \dot{f}_{j}(s)$  is a non-decreasing function in $s$, 
by (2.4), 
 we now conclude 
for all $\{\mu_j\}\in\mathcal{M}$ that
$$
F_1(\{\mu_j\}) \; =\; \lim_{s\to -\infty} \{
\varliminf_{j\to \infty}  \dot{f}_{j}(s)\}
\; 
 \leq \; \varliminf_{j\to \infty} \dot{f}_{j}(0) 
\; =\; 0,
$$
as required. This completes the proof of Theorem.

\bigskip\noindent
{\footnotesize
{\sc Department of Mathematics}\newline
{\sc Osaka University} \newline
{\sc Toyonaka, Osaka, 560-0043}\newline
{\sc Japan}}
\end{document}